\documentclass[10pt]{amsart}

\usepackage{amsmath,amsfonts,amssymb,amsthm}
\usepackage{mathdots}

\newtheorem{theorem}{Theorem}[section]
\newtheorem{lemma}[theorem]{Lemma}
\newtheorem{proposition}[theorem]{Proposition}
\newtheorem{corollary}[theorem]{Corollary}

\theoremstyle{definition}

\theoremstyle{remark}
\newtheorem{remark}[theorem]{Remark}

\numberwithin{equation}{section}

\newcommand{\wo}{\mathbb{E}}
\newcommand{\wwo}[2]{\mathbb{E}\left[{#1}\big|{#2}\right]}
\newcommand{\nat}{\mathbb{N}}
\newcommand{\calk}{\mathbb{Z}}
\newcommand{\rzecz}{\mathbb{R}}
\newcommand{\zesp}{\mathbb{C}}
\newcommand{\saj}[1]{{\mathcal{F}}_{#1}}
\newcommand{\sad}[2]{{\mathcal{F}}_{\le #1,\ge #2}}
\newcommand{\norma}[2]{{\lVert{#2}\rVert}_{#1}}                                 
\DeclareMathOperator{\corr}{corr}                                               
\DeclareMathOperator{\Var}{Var}                                                 
\newcommand{\wVar}[2]{\Var\left[\left.{#1}\right|{#2}\right]}                   

\begin{document}

\title[Generalized random fields]{Generalized stationary random fields with linear regressions - an operator
approach}

\author[W.Matysiak]{Wojciech Matysiak}
\address{Wydzia{\l} Matematyki i Nauk Informacyjnych\\
Politechnika Warszawska\\
Pl. Politechniki 1\\
00-661 Warszawa, Poland\\
and Department of Mathematical Sciences\\
University of Cincinnati\\
Cincinnati, OH 45221-0025\\
USA} \email{matysiak@mini.pw.edu.pl}

\author[P.J.Szab\l owski]{Pawe{\l} J. Szab\l owski}
\address{Wydzia{\l} Matematyki i Nauk Informacyjnych\\
Politechnika Warszawska\\
Pl. Politechniki 1\\
00-661 Warszawa, Poland} \email{pszablowski@elka.pw.edu.pl}
\thanks{}



\date{}

\dedicatory{}

\begin{abstract}
Existence, $L^2$-stationarity and linearity of conditional expectations $\wwo{X_k}{\ldots,X_{k-2},X_{k-1}}$ of square
integrable random sequences $\mathbf{X}=\left( X_{k}\right)_{k\in\mathbb{Z}}$ satisfying
\[
\wwo{X_k}{\ldots,X_{k-2},X_{k-1},X_{k+1},X_{k+2},\ldots}=\sum_{j=1}^\infty b_j\left(X_{k-j}+X_{k+j}\right)
\]
for a real sequence $\left(b_n\right)_{n\in\nat}$, is examined. The analysis is reliant upon the use of Laurent and
Toeplitz operator techniques.
\end{abstract}

\maketitle

\section{Introduction}
Bryc \cite{brycannals} examined square integrable stationary random sequences $\mathbf{X}=\left(
X_{k}\right)_{k\in\mathbb{Z}}$ such that for all integer $k$
\begin{equation}\label{linregr}
\wwo{X_{k}}{\saj{\neq k}}=a\left( X_{k-1}+X_{k+1}\right),
\end{equation}
where $\saj{\neq k}:=\sigma\left(  X_{j}:j\neq k\right)$, $a\in\rzecz$, and
\begin{equation*}
\wwo{X_k^2}{\mathcal{F}_{\neq k}}=Q(X_{k-1},X_{k+1})
\end{equation*}
for some (symmetric) quadratic form $Q$. The sequences were further analyzed in \cite{Bryc99}, \cite{matszab2},
\cite{matszab3} and almost complete characterization of their finite dimensional distributions is known. Bryc's random
fields are related to Hammersley's harnesses \cite{hammersley}, recently revived by Mansuy and Yor \cite{mansuyyor},
and Bryc and Weso{\l}owski \cite{brycwesolo1} (see also \cite{bmw1} and references therein).

The first step of the analysis in \cite{brycannals} that led to knowing the distributional structure of $\mathbf{X}$
was to observe that \eqref{linregr} implies that one-sided regressions are linear:
\begin{equation}\label{onesidedannals}
\wwo{X_k}{\saj{k-1}}=\alpha X_{k-1},
\end{equation}
where $\saj{k-1}:=\sigma\left(  X_{j}:j\le k-1\right)$ and $\alpha\in\rzecz$ ($\alpha$ is easily identified as the
correlation coefficient $\corr(X_{k},X_{k-1})$).

This paper is devoted to the introductory analysis of random sequences that satisfy a generalization of condition
\eqref{linregr}. Namely, throughout the paper $\mathbf{X}$ will denote square integrable and normalized random sequence
satisfying for all $k\in\calk$
\begin{equation}\label{infinitelr}
\wwo{X_k}{\saj{\ne k}}=\sum_{j=1}^\infty b_j\left(X_{k-j}+X_{k+j}\right)
\end{equation}
for some real sequence $\left(b_n\right)_{n\in\nat}$ (here and further it is assumed that all the series involving
elements of $\mathbf{X}$ converge in $L^2$). We shall formulate conditions for sequences $\left(b_n\right)_n$ under
which processes defined by \eqref{infinitelr} exist, are $L^2$-stationary and have linear one-sided regressions:
\begin{equation}
\wwo{X_k}{\saj{\le k-1}}=\sum_{j=1}^\infty \beta_j X_{k-j}
\end{equation}
for some coefficients $\left(\beta_n\right)_n$ related to $(b_n)_n$. (It will turn out that if there exists $N\in\nat$
such that $b_j=0$ for $j>N$, then $\beta_j=0$ for $j>N$.)

Similar problems were examined in a much greater generality by Williams \cite{williams1}, followed by Kingman
\cite{kingman1}, \cite{kingman2}. In the case of generalized Bryc's random fields, the matrix of linear regression
coefficients, denoted later as $L(b)$, is a (symmetric) Laurent matrix, i.e. doubly infinite matrix, which is constant
along the diagonals. The structure of the matrix $L(b)$ allows stronger results for the existence of $\mathbf{X}$ to be
obtained than the ones obtained in \cite{williams1} and \cite{kingman2}. In our analyses, we will also refer to
Toeplitz matrices, which are infinite (but not doubly infinite) matrices with constant diagonals.

Section \ref{sec:mainresults} is the central section of the paper. Main results of the paper are listed in subsections
\ref{subsec:existence} and \ref{subsec:linearity}, while subsections \ref{subsec:Laurent} and \ref{subsec:Toeplitz}
gather the important facts from the theory of Laurent and Toeplitz operators needed in further considerations. Section
\ref{sec:proofs} contains proofs and auxiliary results. Concluding remarks can be found in Section
\ref{sec:concluding}.

\section{Main results}\label{sec:mainresults}

\subsection{Laurent matrices and their symbols}\label{subsec:Laurent}
Given a sequence of complex numbers $(a_n)_{n\in\calk}$, one can construct a Laurent matrix
$$
\left[a_{m-n}\right]_{m,n\in\mathbb{Z}}=
\begin{pmatrix}
\ddots & \ddots & \ddots & \ddots & \ddots \\
\ddots & a_0 & a_{-1} & a_{-2} & \ddots \\
\ddots & a_1 & \framebox[1.2\width][c]{$a_0$} & a_{-1} & \ddots \\
\ddots & a_2 & a_1 & a_0 & \ddots \\
\ddots & \ddots & \ddots & \ddots & \ddots
\end{pmatrix}.
$$
(Throughout the paper, the box indicates the entry in $(0,0)$ position in the case of doubly infinite matrices, or $0$
position in the case of doubly infinite vectors.) By the classical theorem of Toeplitz \cite{Toeplitz}, the matrix
defines a bounded operator on $l^2(\calk)$ if and only if the numbers $(a_n)_n$ are the Fourier coefficients of some
function $a\in L^\infty(\mathbb{T})$ (where $\mathbb{T}$ denotes the complex unit circle)
\[
a_n=\frac{1}{2\pi}\int_{-\pi}^{\pi} a\left(e^{i\theta}\right)e^{-i n \theta} d\theta, \quad n\in\calk.
\]
If such a function exists then it is unique and called the symbol of the Laurent matrix $\left[a_{m-n}\right] _{m,n}$
(for a readable introduction to the theory of Laurent and Toeplitz operators on $l^p$ spaces, see
\cite{BottcherSilbermann}). We shall denote by $L(a)$ both Laurent matrix $\left[a_{m-n}\right] _{m,n}$ and the bounded
operator generated by it.

It is known that if $a\in L^\infty(\mathbb{T})$, then the spectrum of the operator $L(a)$ is equal to the spectrum
$\mathcal{R}(a)$ of $a$ as an element of Banach algebra $L^\infty(\mathbb{T})$, which in turn is equal to the essential
range of $a$ ($|\cdot|$ stands for the Lebesgue measure):
\begin{equation}\label{essrange}
\mathcal{R}(a)=\bigl\{\lambda\in\zesp:\bigl|\left\{t\in\mathbb{T}:|a(t)-\lambda|<\epsilon\right\}\bigr|>0\ \forall
\epsilon>0\bigr\},
\end{equation}
and if $0\notin \mathcal{R}(a)$, then the inverse of $L(a)$ is the Laurent matrix with symbol $a^{-1}$ (Theorem 1.2,
\cite{BottcherSilbermann}).

\subsection{Existence and $L^2$-stationarity}\label{subsec:existence}

For a given real sequence $(b_n)_{n\in\nat}$, consider a doubly infinite matrix with $(i,j)$-th entry ($i,j\in\calk$)
defined as
\begin{equation}\label{basicmatrix}
\left\{%
\begin{array}{ll}
    -b_{|i-j|}, & \hbox{if $i\neq j$;} \\
    1, & \hbox{if $i=j$.} \\
\end{array}%
\right.
\end{equation}

We will denote the symbol of the matrix defined in \eqref{basicmatrix} (if it exists) by $b$, and the matrix itself as
well as the corresponding operator by $L(b)$. 

\begin{proposition}\label{existence}
If $\left(b_n\right)_{n\in\nat}$ is a sequence of real numbers such that the symbol $b\in L^\infty(\mathbb{T})$ is
positive and $0\notin \mathcal{R}(b)$, then there exist square integrable random sequences
$\mathbf{X}=\left(X_k\right)_{k\in\calk}$ satisfying \eqref{infinitelr} for each $k\in\calk$.
\end{proposition}
Since the essential range of continuous function $a$ is the image $a(\mathbb{T})$, we see that the sufficient condition
for the existence of $\mathbf{X}$ in the case of $b\in C(\mathbb{T})$ is $b>0$.

\begin{remark}\label{brycannals}
If $b_1=a$ and $b_n=0$ for $n\ge2$, one obtains the case considered by Bryc \cite{brycannals}. It is easy to verify
that $b(\theta)=1-2a\cos\theta$ is positive for $\theta\in (-\pi,\pi]$ if and only if $|a|<\frac{1}{2}$. Since $b$ is
continuous, we arrive at the assertion of Theorem 1 \cite{matszab2} (see also \cite{brycannals}).
\end{remark}

\begin{remark}\label{williamsremark1}
Williams \cite{williams1} considered, among other things, the problem of the existence of square integrable random
sequences $\left(X_k\right)_k$, with $k$ belonging to an arbitrary countable set and satisfying
\begin{equation}\label{williamstypeequation}
\wwo{X_k}{\sigma\{X_j:j\ne k\}}=\sum_{j\ne k} a_{k,j}X_j,
\end{equation}
with $a_{k,j}\ge0$, $\sum_{j} a_{k,j}\le 1$ and the sets $\{j:a_{k,j}>0\}$ being finite for all $k$. In the case of
symmetry of the matrix $[a_{k,j}]_{k,j}$, Williams proved that the sequences exist if the matrix $[a_{k,j}]_{k,j}$ is
invertible and positive definite. Thus Proposition \ref{existence} extends Williams's result to the case of an infinite
number of non-zero coefficients in \eqref{williamstypeequation}, with the price paid of assuming the special structure
of the matrix $[a_{k,j}]_{k,j}$.
\end{remark}

\begin{remark}
Kingman \cite{kingman1} considered the problem of the existence of finite sets of $L^1$ random variables satisfying
\eqref{williamstypeequation}; in \cite{kingman2} the problem was extended to random sequences under the assumption that
the matrix $[|a_{k,j}|]_{k,j}$ is transient.
\end{remark}

\begin{proposition}\label{L2stationarity}
If $\mathbf{X}$ satisfies \eqref{infinitelr} for a sequence $\left(b_n\right)_{n}$ such that $b\in
L^\infty(\mathbb{T})$ is positive and $0\notin \mathcal{R}(b)$, then $\mathbf{X}$ is $L^2$-stationary.
\end{proposition}

From now on, we will denote the correlation coefficients of $\mathbf{X}$ by $r_{|k|}=\wo X_0 X_k$ for $k\in\calk$.

\subsection{Toeplitz matrices with symbols in Wiener algebra}\label{subsec:Toeplitz}
Throughout this subsection, we will consider matrices $L(b)$ with symbols in Wiener algebra (see
\cite{BottcherSilbermann}), i.e. in the set $W=W(\mathbb{T})$ of all complex-valued functions $a$ such that
\[
a(t)=\sum_{k\in\calk} a_k t^k,\ \textrm{where}\ \sum_{k\in\calk} |a_k|<\infty\ \textrm{and} \
t=\exp(i\theta)\in\mathbb{T}.
\]
$W$ is a Banach algebra with pointwise algebraic operations and the norm $\norma{}{a}=\sum_{k\in\calk} |a_k|$. Clearly,
$W(\mathbb{T})\subset C(\mathbb{T})$, so if $b\in W(\mathbb{T})$ is positive, the corresponding random sequence
$\mathbf{X}$ exists. An important result concerning Wiener algebra is Wiener's theorem: if $a\in W(\mathbb{T})$ and $a$
has no zeros on $\mathbb{T}$ then $a^{-1}=1/a\in W(\mathbb{T})$. Denoting by $\mathcal{G}(A)$ the set of all invertible
elements of a Banach algebra $A$, we can rephrase Wiener's theorem as
\begin{equation}\label{wiener}
\mathcal{G}\left(W(\mathbb{T})\right)=\left\{a\in W: a(t)\ne 0\ \forall t\in\mathbb{T}\right\}.
\end{equation}

The reason for considering absolutely summable sequences $\left(b_n\right)_n$ is that in the proof of Theorem
\ref{onesidedthm}, one needs to use the Toeplitz operator $T(a)$ generated by an absolutely summable sequence
$\left(a_k\right)_{k\in\calk}$.

The Toeplitz matrix (see \cite{BottcherSilbermann}) defined by a sequence $\left(c_k\right)_{k\in\calk}$ of complex
numbers is the infinite matrix
$$
[c_{m-n}]_{m,n\in\nat}=
\begin{pmatrix}
c_0 & c_{-1} & c_{-2} & \ldots \\
c_1 & c_0 & c_{-1} & \ldots \\
c_2 & c_1 & c_0 & \ldots \\
\vdots & \ddots & \ddots & \ddots
\end{pmatrix}
.
$$
It is known that if $\sum_{k\in\calk} |c_k|<\infty$, then $[c_{m-n}]_{m,n}$ induces a bounded operator on $l^1(\nat)$
(Proposition 7.1, \cite{BottcherSilbermann}). Analogously to Laurent operators, function
\[
c(t)=\sum_{k\in\calk} c_k t^k, \quad t=\exp(i\theta)\in\mathbb{T}
\]
is called the symbol of the Toeplitz operator. Also, analogously to the Laurent case, we shall denote both the matrix
and the operator corresponding to $c$ as $T(c)$. By theorems of Gohberg and Duduchava (Theorems 7.3 and 7.4,
\cite{BottcherSilbermann}), if $c\in W(\mathbb{T})$, then $T(c)$ is invertible on $l^1(\nat)$ if and only if $0\notin
c(\mathbb{T})$ and $\textrm{wind}(c,0)=0$, where $\textrm{wind}(c,0)$ is the winding number of $c$ with respect to the
origin.

\subsection{Linearity of one-sided regressions}\label{subsec:linearity}
\begin{theorem}\label{onesidedthm}
If $\mathbf{X}$ satisfies \eqref{infinitelr} and $0<b\in W(\mathbb{T})$, then one-sided regressions are linear:
\begin{equation}\label{onesidedeqn}
\wwo{X_k}{\saj{\le{k-1}}}=\sum_{j=1}^\infty \beta_j X_{k-j}.
\end{equation}
Coefficients $\left(\beta_n\right)_{n\in\nat}$ satisfy the following equation: if
\[
\beta(t)=1-\sum_{j=1}^\infty \beta_j t^j,\quad \widetilde{\beta}(t)=\beta(1/t),\quad t\in\mathbb{T},
\]
and
\begin{equation}
v=1-2\sum_{j=1}^\infty b_j r_j,\quad w=1-\sum_{j=1}^\infty \beta_j r_j,
\end{equation} then
\begin{equation}\label{bbetasymbols}
b(t)=\frac{v}{w}\beta(t)\widetilde{\beta}(t),\quad t\in\mathbb{T}.
\end{equation}
\end{theorem}
\begin{remark}
Let $H^\infty$ ($\overline{H^\infty}$) denote the closed subalgebra of the Banach algebra $L^\infty(\mathbb{T})$
consisting of all functions with vanishing Fourier coefficients with negative (positive) indices. It is well known that
in general (unlike to the Laurent case), the product of two Toeplitz operators is not a Toeplitz operator.
Nevertheless, if only $a_1\in\overline{H^\infty}$ (so $T(a_1)$ is upper triangular), $a_2\in L^\infty$ and $a_3\in
H^\infty$ (so $T(a_3)$ is lower triangular), then
\begin{equation}\label{toeplitzproduct}
T\left(a_1\right)T\left(a_2\right)T\left(a_3\right)=T\left(a_1a_2a_3\right)
\end{equation}
(see Proposition 1.13 \cite{BottcherSilbermann}). In view of the above, \eqref{bbetasymbols} can be interpreted as
\[
\begin{pmatrix}
1 & -b_{1} & -b_{2} & \ldots \\
-b_1 & 1 & -b_{1} & \ldots \\
-b_2 & -b_1 & 1 & \ldots \\
\vdots & \ddots & \ddots & \ddots
\end{pmatrix}
=\frac{v}{w}
\begin{pmatrix}
1 & -\beta_{1} & -\beta_{2} & \ldots \\
0 & 1 & -\beta_{1} & \ldots \\
0 & 0 & 1 & \ldots \\
\vdots & \ddots & \ddots & \ddots
\end{pmatrix}
\begin{pmatrix}
1 & 0 & 0 & \ldots \\
-\beta_1 & 1 & 0 & \ldots \\
-\beta_2 & -\beta_1 & 1 & \ldots \\
\vdots & \ddots & \ddots & \ddots
\end{pmatrix}
.
\]
\end{remark}

\begin{remark}
Note that $v=\wVar{X_k}{\saj{\ne k}}$ and $w=\wVar{X_k}{\saj{\le k-1}}$.
\end{remark}

In the finite case, that means with only a finite number of non-zero elements of $\left(b_n\right)_n$, it is possible
to give a more explicit description of the connection between coefficients $\left(b_n\right)_n$ and
$\left(\beta_n\right)_n$. (In such case, the symbol $b$ is a trigonometric polynomial and $L(b)$ is a band Laurent
operator.)

\begin{corollary}\label{onesidedfinite}
If $\mathbf{X}$ satisfies \eqref{infinitelr} and $N:=\sup\{n\in\nat: b_n\ne0\}<\infty$, then \eqref{onesidedeqn} holds
and $\sup\{n\in\nat: \beta_n\ne0\}=N$. Furthermore,
\begin{eqnarray}\label{beta2b}
b_1 &=& \frac{\beta_{N-1}-\beta_1\beta_N}{1+\beta_1^2+\ldots+\beta_N^2},\nonumber\\
b_2 &=& \frac{\beta_{N-2}-\beta_1\beta_{N-1}-\beta_2\beta_{N-2}}{1+\beta_1^2+\ldots+\beta_N^2},\nonumber\\
\vdots & & \vdots\\
b_{N-1} &=& \frac{\beta_1-\beta_1\beta_2-\ldots-\beta_{N-1}\beta_N}{1+\beta_1^2+\ldots+\beta_N^2},\nonumber\\
b_N &=& \frac{\beta_N}{1+\beta_1^2+\ldots+\beta_N^2}\nonumber.
\end{eqnarray}
\end{corollary}

\begin{remark}
If $\mathbf{X}$ satisfies \eqref{infinitelr} with $b_1=a$ and $b_n=0$ for $n\ge2$ (see \cite{brycannals}), then by
Corollary \ref{onesidedfinite} one gets that $\wwo{X_k}{\saj{\le k-1}}=\beta_1 X_{k-1}$ (cf. \eqref{onesidedannals}).
Multiplying this equation by $X_{k-1}$ and taking the expectations, yields $\beta_1=r_1$. Finally, $a=r_1/(1+r_1^2)$ by
\eqref{beta2b} (see Theorem 3.1 \cite{brycannals}, and paper \cite{matszab2}).
\end{remark}

\begin{remark}
In \cite{williams1}, Williams proved that under the assumptions stated in Remark \ref{williamsremark1}, a variant of
\eqref{williamstypeequation} holds with the $\sigma$-algebra $\sigma\{X_j:j\ne k\}$ being replaced by more general
$\sigma$-algebras $\sigma\{X_j:j\notin S\}$, $S\subset J$. Thus, Theorem \ref{onesidedthm} is an extension of
Williams's result to the case of $L^2$-stationary random sequences indexed by integers and with an infinite number of
non-zero regression coefficients.
\end{remark}

\section{Proofs of main results}\label{sec:proofs}

\subsection{Existence and $L^2$-stationarity}
The main observation behind the proof of Proposition \ref{existence} is the following relation of the matrix of linear
regression coefficients $L(b)$ and the covariance matrix of the random sequence $\mathbf{X}$ that can be traced (in the
finite case) to at least the paper of Kingman \cite{kingman1}; in implicit form, one can find it also in
\cite{williams1}. Namely, if $\mathbf{X}$ satisfies \eqref{infinitelr} and $\mathbf{R}$ denotes its (bi-infinite)
covariance matrix, then
\begin{equation}\label{basicrelation}
L(b)\ \mathbf{R}=\mathbf{\Delta},
\end{equation}
where $\mathbf{\Delta}$ is the (bi-infinite) diagonal matrix with $(j,j)$-th entry equal to
\begin{equation}\label{definitionv}
v_j=\wo\left[X_j-\wwo{X_j}{\saj{\ne j}}\right]^2,\quad j\in\calk.
\end{equation}
It is easy to check that \eqref{basicrelation} holds true by multiplying \eqref{infinitelr} by $X_j$, $j\in\calk$, and
calculating the unconditional expectation (see \cite{matszab2}, Proposition 2).

\begin{proof}[Proof of Proposition \ref{existence}]
Positivity of $b$ implies positive definiteness of $L(b)$. Condition $0\notin \mathcal{R}(b)$ guarantees existence of
the inverse $L^{-1}(b)=L(b^{-1})$, which is symmetric and positive definite. Therefore, there exists the centered
Gaussian sequence $\mathbf{X}=\left(X_k\right)_{k\in\calk}$ with the covariance matrix $L^{-1}(b)$. $L^{-1}(b)$ is a
Laurent matrix, hence $\mathbf{X}$ is $L^2$-stationary. Therefore, all the numbers $v_j$ defined in \eqref{definitionv}
are equal to a real number $v$. Since Gaussian sequences have linear regressions, $\mathbf{X}$ satisfies
\eqref{infinitelr}.
\end{proof}
\begin{proof}[Proof of Proposition \ref{L2stationarity}]
Matrix $L^{-1}(b)$ in the proof of Proposition \ref{existence} is Laurent.
\end{proof}
In view of Corollary \ref{L2stationarity}, we can write the covariance matrix as $\mathbf{R}=L(r)$ with $r\in
L^{\infty}(\mathbb{T})$, and then \eqref{basicrelation} reads as $L(b)L(r)=v\mathbf{I}$ (here and further, $\mathbf{I}$
denotes the identity matrix of appropriate dimensions). Since the symbol of the product of two Laurent matrices is the
product of the symbols, we get that
\begin{equation}\label{basicrelationsymbols}
b(t)r(t)=v,\  t=\exp(i\theta)\in\mathbb{T}.
\end{equation}

\subsection{Linearity of one-sided regressions}
Lemma \ref{easyaf} is needed in the proof of Theorem \ref{onesidedthm} (the theorem imitates the ideas of Williams
\cite{williams1}, for a simple case see \cite{williams2}, Chapter 15.10).

\begin{lemma}\label{easyaf}
Let $\left(X,\norma{}{\cdot}\right)$ be a Banach space and $\mathcal{B}(X)$ be the Banach algebra of bounded linear
operators on $X$ (the operator norm will also be denoted by $\norma{}{\cdot}$).

Suppose $\left(y_n\right)_{n\in\nat}\in X$ and $\left(T_n\right)_{n\in\nat}\in \mathcal{B}(X)$. Assume that $y_n\to
y_\infty$ in $X$, $T_n\to T_\infty$ in $\mathcal{B}(X)$ and $T_\infty$ is invertible. If $\left(x_n\right)_{n\in\nat}$
satisfies $T_n x_n=y_n$, then $x_n\to x_\infty$ in $X$, where $x_\infty=T^{-1}_\infty y_\infty$.
\end{lemma}

\begin{proof}
Since for any (non-zero) Banach algebra $A$ with identity, $\mathcal{G}(A)$ is open in $A$ (see e.g. \cite{rudin},
18.4, Corollary 1), then there exists $N\in\nat$ such that $T_n^{-1}\in\mathcal{B}(X)$ for all $n>N$. For those $n$,
\begin{multline*}
\norma{}{x_n-x_\infty}=\norma{}{T_n^{-1}y_n -T_\infty^{-1} y_\infty}=\norma{}{T_n^{-1}y_n -T_\infty^{-1}y_n +
T_\infty^{-1} y_n - T_\infty^{-1} y_\infty}\le \\
\le\norma{}{T_n^{-1}-T_\infty^{-1}}\norma{}{y_n}+\norma{}{T_\infty^{-1}}\norma{}{y_n-y_\infty}.
\end{multline*}
For any Banach algebra $A$, the mapping $x\mapsto x^{-1}$ is a homeomorphism of $\mathcal{G}(A)$ onto itself
(\cite{rudin},18.4, Corollary 1), so $\norma{}{T_n^{-1}-T_\infty^{-1}}\to0$. Clearly, due to the convergence of
$\left(y_n\right)_n$, the sequence $\left(\norma{}{y_n}\right)_n$ is bounded, so the right hand side of the above
inequality converges to $0$.
\end{proof}

\begin{proof}[Proof of Theorem \ref{onesidedthm}]
Fix $k\in\calk$, $m\in\nat$ and set $y_n:=\wwo{X_{k+n-1}}{\sad{k-1}{k+m}}$, $n=1,\ldots,m$. Since
\begin{multline*}
y_n=\wwo{\wwo{X_{k+n-1}}{\saj{\ne k+n-1}}}{\sad{k-1}{k+m}}=\\
=\sum_{j=n}^\infty b_j X_{k+n-1-j} +\sum_{j=1}^{n-1} b_j \wwo{X_{k+n-1-j}}{\sad{k-1}{k+m}}+\\  +\sum_{j=1}^{m-n} b_j
\wwo{X_{k+n-1+j}}{\sad{k-1}{k+m}} + \sum_{j=m-n+1}^\infty b_j X_{k+n-1+j}
\end{multline*}
for $n=1,\ldots,m$, we see that the vector $\mathbf{y}=\begin{pmatrix}y_1 & \ldots & y_m\end{pmatrix}^T$ satisfies the
system of linear equations $T_m(b)\ \mathbf{y}=\mathbf{d}$, with the main matrix
\[
T_m(b)=
\begin{pmatrix}
1 & -b_1 & -b_2 & \ldots & -b_{m-1}\\
-b_1 & 1 & -b_1 & \ldots & -b_{m-2}\\
-b_2 & -b_1 & 1 & \ldots & -b_{m-3}\\
\ldots & \ldots & \ldots & \ldots & \ldots\\
-b_{m-1} & -b_{m-2} & -b_{m-3} & \ldots & 1
\end{pmatrix},
\]
and vector $\mathbf{d}$ consisting of some infinite linear combinations (convergent in $L^2$) of random variables
$\ldots,X_{k-2},X_{k-1}$ and $X_{k+m}$, $X_{k+m+1},\ldots$, with some real coefficients that do not depend on $k$.

We claim that (for sufficiently large $m$), $T_m(b)$ is invertible. Indeed, $T_m(b)$ is the principal $m\times m$
submatrix of the infinite Toeplitz matrix $T(b)$. The symbol $b$ is real-valued, continuous and positive (hence
$\textrm{wind}(b,0)=0$), so by the Krein theorem (\cite{BottcherSilbermann}, Theorem 1.15), $T(b)$ is invertible on
$l^2$. Hence, and by the assumption that $b\in C(\mathbb{T})$, we are in a position to use the Gohberg-Feldman theorem
(\cite{BottcherSilbermann}, Theorem 2.11). This theorem implies that the sequence of $n\times n$ matrices
$\left(T_n(b)\right)_{n\in\nat}$ is stable, which means, among other things, that $T_n(b)$ is invertible for
sufficiently large $n$.

Therefore, solving the system $T_m(b)\ \mathbf{y}=\mathbf{d}$, we get, in particular,
\begin{equation}\label{indukcjaeqn}
y_1=\wwo{X_k}{\sad{k-1}{k+m}}=\sum_{j=1}^\infty \beta_{j,m} X_{k-j}+ \sum_{j=1}^\infty \gamma_{j,m} X_{k+m-1+j},
\end{equation}
say, with $\beta_{j,m}$ and $\gamma_{j,m}$ real and not depending on $k$. Multiplying both of sides of
\eqref{indukcjaeqn} by $\ldots,X_{k-2},X_{k-1}$ and $X_{k+m},X_{k+m+1},\ldots$, and taking expectations, we obtain an
infinite system of linear equations with unknowns $\beta_{j,m}$ and $\gamma_{j,m}$ ($j\in\nat$), which can be written
in the matrix form $\mathbf{A}(m) \mathbf{x}(m)= \mathbf{b}(m)$ with
\begin{eqnarray*}
\mathbf{b}(m)&=&{\begin{pmatrix}\ldots & r_2 & \framebox[1.2\width][c]{$r_1$} & r_m & r_{m+1} & \ldots\end{pmatrix}}^T,\\
\mathbf{x}(m)&=&{\begin{pmatrix}\ldots & \beta_{2,m} & \framebox[1.1\width][c]{$\beta_{1,m}$} & \gamma_{1,m} & \gamma_{2,m} & \ldots\end{pmatrix}}^T,\\
\mathbf{A}(m)&=&\begin{pmatrix}
\mathbf{A}_{1} & \mathbf{A}_{2}(m) \\
\mathbf{A}_{3}(m) & \mathbf{A}_{4}
\end{pmatrix},
\end{eqnarray*}
and
\begin{eqnarray*}
\mathbf{A}_{1}\phantom{(m)}=
\begin{pmatrix}
  \ddots & \vdots & \vdots & \vdots \\
  \cdots & 1 & r_1 & r_2 \\
  \cdots & r_1 & 1 & r_1 \\
  \cdots & r_2 & r_1 & \framebox[1.2\height][c]{$1$} \\
\end{pmatrix}
&,& \mathbf{A}_{2}(m)=
\begin{pmatrix}
  \vdots & \vdots & \iddots \\
  r_{m+3} & r_{m+2} & \cdots \\
  r_{m+2} & r_{m+3} & \cdots \\
  r_{m+1} & r_{m+2} & \cdots \\
\end{pmatrix}
,\\
\mathbf{A}_{3}(m)=
\begin{pmatrix}
  \cdots & r_{m+2} & r_{m+1} \\
  \cdots & r_{m+3} & r_{m+2} \\
  \cdots & r_{m+4} & r_{m+3} \\
  \iddots & \vdots & \vdots \\
\end{pmatrix}
&,& \mathbf{A}_{4}\phantom{(m)}=
\begin{pmatrix}
  1 & r_1 & r_2 & \cdots \\
  r_1 & 1 & r_1 & \cdots \\
  r_2 & r_1 & 1 & \cdots \\
  \vdots & \vdots & \vdots & \ddots \\
\end{pmatrix}.
\end{eqnarray*}
Observe that submatrices $\mathbf{A}_{1}$ and $\mathbf{A}_{4}$ do not depend on $m$.

By \eqref{basicrelationsymbols} and by the Wiener theorem, $r\in W(\mathbb{T})$. Therefore,
\[
\mathbf{b}(m)\longrightarrow{\begin{pmatrix}\ldots & r_2 & \framebox[1.2\width][c]{$r_1$} & 0 & 0 &
\ldots\end{pmatrix}}^T=:\mathbf{b}_\infty\quad \textrm{in}\ l^1(\calk).
\]

Since $ \norma{}{\mathbf{A}}=\sup_{n\in\calk} \sum_{m\in\calk} \left|a_{m,n}\right| $ for
$\mathbf{A}=\left[a_{m,n}\right]_{m,n\in\calk}$, regarded as an element of $\mathcal{B}\left(l^1(\calk)\right)$,
$\mathbf{A}_m$ tends in $\mathcal{B}\left(l^1(\calk)\right)$ to $\mathbf{A}_\infty$ defined as
\[
\mathbf{A}_\infty=\begin{pmatrix}
  \mathbf{A}_{1} & \mathbf{0} \\
  \mathbf{0} & \mathbf{A}_{4} \\
\end{pmatrix}
\]
(here and further, $\mathbf{0}$ denotes the zero matrix of the appropriate infinite dimensions).

Now we will show that the inverse of $\mathbf{A}_\infty$ exists and belongs to $\mathcal{B}\left(l^1(\calk)\right)$.
First, observe that $\mathbf{A}_{4}$ can be identified in a natural way with an operator acting on $l^1(\nat)$, namely
with the Toeplitz operator generated by sequence $\left(a_k\right)_{k\in\calk}$ with $a_{k}=r_{|k|}$ for $k\in\calk$,
having the symbol $r$ equal to the symbol of the correlation matrix $L(r)$. By \eqref{basicrelationsymbols}, $r$ is
real-valued, continuous and positive (hence $\textrm{wind}(r,0)=0$), so by the Gohberg and Duduchava theorems,
$\mathbf{A}_{4}$ is invertible on $l^1$. Analogously, $\mathbf{A}_{1}$ is invertible on $l^1$. Clearly,
$\mathbf{A}^{-1}_\infty\in \mathcal{B}\left(l^1(\calk)\right)$ and is of the form
\[
\mathbf{A}^{-1}_\infty=\begin{pmatrix}
  \mathbf{A}_1^{-1} & \mathbf{0} \\
  \mathbf{0} & \mathbf{A}_4^{-1} \\
\end{pmatrix}.
\]

Now we are in a position to use Lemma \ref{easyaf} with $X=l^1(\calk)$. Doing so, we get that
\[
\mathbf{x}(m)\longrightarrow \mathbf{A}^{-1}_\infty \mathbf{b}_\infty=:\mathbf{x}_\infty \quad \textrm{in}\ l^1(\calk),
\]
and
\[
\mathbf{x}_\infty={\begin{pmatrix}\ldots & \beta_{2} & \framebox[1.2\width][c]{$\beta_{1}$} & 0 & 0 &
\ldots\end{pmatrix}}^T
\]
for some real $\beta_{j}$, $j\in\nat$. Note that
\begin{equation}\label{betas}
\begin{pmatrix}
\beta_1\\
\beta_2\\
\beta_3\\
\vdots
\end{pmatrix}
=
\begin{pmatrix}
  1 & r_1 & r_2 & r_3 &\cdots \\
  r_1 & 1 & r_1 & r_2 &\cdots \\
  r_2 & r_1 & 1 & r_1 & \cdots \\
  \vdots & \vdots & \vdots & \vdots & \ddots \\
\end{pmatrix}^{-1}
\begin{pmatrix}
r_1 \\
r_2 \\
r_3 \\
\vdots
\end{pmatrix}.
\end{equation}

On the other hand, by applying L\'{e}vy's convergence theorem for backward martingales to \eqref{indukcjaeqn}, we
deduce that the limit
\[
L:=\lim_{m\to\infty} \Bigl( \sum_{j=1}^\infty \beta_{j,m} X_{k-j}+\sum_{j=1}^\infty \gamma_{j,m} X_{k+m+j-1} \Bigr)
\]
exists almost surely and in $L^1$. We claim that $L=\sum_{j=1}^\infty \beta_{j} X_{k-j}$. Indeed,
\begin{multline*}
\wo\left|\sum_{j=1}^\infty \beta_{j,m} X_{k-j}+\sum_{j=1}^\infty \gamma_{j,m}
X_{k+m+j-1}-\sum_{j=1}^\infty \beta_{j} X_{k-j}\right|\le\\
\le \wo\left|\sum_{j=1}^\infty \beta_{j,m} X_{k-j}-\sum_{j=1}^\infty \beta_{j}
X_{k-j}\right|+\wo\left|\sum_{j=1}^\infty \gamma_{j,m} X_{k+m+j-1}\right|.
\end{multline*}
But
\[
\wo\left|\sum_{j=1}^\infty \gamma_{j,m} X_{k+m+j-1}\right|\le\sum_{j=1}^\infty
\left|\gamma_{j,m}\right|\wo\left|X_{k+m-j-1}\right|\le\sum_{j=1}^\infty \left|\gamma_{j,m}\right|\longrightarrow 0
\]
and
\begin{multline*}
\wo\left|\sum_{j=1}^\infty \beta_{j,m} X_{k-j}-\sum_{j=1}^\infty \beta_{j,\infty} X_{k-j}\right|\le \sum_{j=1}^\infty
\left|\beta_{j,m}-\beta_{j,\infty}\right|\wo\left|X_{k-j}\right|\le\\
\le \sum_{j=1}^\infty \left|\beta_{j,m}-\beta_{j,\infty}\right|\longrightarrow 0,
\end{multline*}
when $m$ tends to infinity, which establishes our claim.

Thus, by L\'{e}vy's theorem applied to \eqref{indukcjaeqn},
\begin{equation*}
\wo\Biggl[{X_k}\left|{\bigcap_{m=1}^\infty \sad{k-1}{k+m}}\right. \Biggr]=\sum_{j=1}^\infty \beta_{j} X_{k-j}.
\end{equation*}
Taking the conditional expectation of both sides of the above equation with respect to $\saj{\le k-1}$, yields
\eqref{onesidedeqn}.

 Since $\widetilde{\beta}\in\overline{H^\infty}$,
$T\left(\widetilde{\beta}\right)T(r)$ is a Toeplitz matrix by \eqref{toeplitzproduct}. Using
\begin{equation}\label{yulewalker}
T(r)
\begin{pmatrix}
\beta_1\\
\beta_2\\
\beta_3\\
\vdots
\end{pmatrix}
=
\begin{pmatrix}
r_1 \\
r_2 \\
r_3 \\
\vdots
\end{pmatrix},
\end{equation}
which is a direct consequence of \eqref{betas}, one can easily verify that the only non-zero element of the first row
of the matrix $T\left(\widetilde{\beta}\right)T(r)$ is its first element, which is equal to $w$. In other words,
$\widetilde{\beta}\ r\in H^\infty$. Since $\beta\in H^\infty$, we get that $\widetilde{\beta}\ r\ \beta\in H^\infty$,
which means, again due to \eqref{toeplitzproduct}, that $T\left(\widetilde{\beta}\right)T(r)T\left(\beta\right)$ is a
lower triangular Toeplitz matrix. But $T\left(\widetilde{\beta}\right)T(r)T\left(\beta\right)$ is symmetric, hence it
must be a diagonal matrix with constant diagonal $w$. Thus, $\beta(1/t)r(t)\beta(t)=w$ and multiplying it by $b$ and
using \eqref{basicrelationsymbols}, we arrive at \eqref{bbetasymbols}.
\end{proof}

\begin{proof}[Proof of Corollary \ref{onesidedfinite}]
If $\sup\{n\in\nat: b_n\ne0\}=N$, then from comparing the powers of $t$ in \eqref{bbetasymbols}, it follows that for
$n\ge N+1$
\[
\beta_n-\sum_{j=1}^\infty \beta_j \beta_{n+j}=0,
\]
which means that the vector $\begin{pmatrix}\beta_{N+1} & \beta_{N+2} & \ldots \end{pmatrix}^T$ belongs to the kernel
of $T(\beta)$. Since $T(\beta)$ is invertible, we get that $\sup\{n\in\nat: \beta_n\ne0\}\le N$. Comparing the
coefficients at $t^N$ at both sides of \eqref{bbetasymbols}, we deduce that $\beta_N\ne0$ and $v/w=b_N/\beta_N$;
comparing all the powers of $t$ in \eqref{bbetasymbols}, we get \eqref{beta2b}.
\end{proof}

\section{Concluding remarks}\label{sec:concluding}

\begin{remark}
Random sequences satisfying \eqref{infinitelr} are the examples of AR time series; equation \eqref{yulewalker} is the
Yule-Walker equation.
\end{remark}

\begin{remark}
Since one of the assumptions of Theorem \ref{onesidedthm} is the absolute summability of the sequence
$\left(b_n\right)_n$, which in view of Wiener's theorem means the absolute summability of the sequence of correlation
coefficients of $\mathbf{X}$, one may rephrase the assertion of Theorem \ref{onesidedthm} by saying that the sequences
with short memory, satisfying \eqref{infinitelr}, have linear one-sided regressions. The question of extending Theorem
\ref{onesidedthm} to the case of processes with long memory is left open.
\end{remark}

\begin{remark}[Conditions for positivity of the symbol $b$] If $N=\sup\{n\in\nat: b_n\ne0\}<\infty$, then a
lemma of Fej\'er and Riesz asserts that cosine polynomial $ b(\theta)=1-2\sum_{n=1}^N b_n \cos (n \theta) $ is
non-negative if and only if it can be expressed in the form
\[
b(\theta)=\left|\sum_{n=0}^N c_n \exp(i n\theta) \right|^2,
\]
where
\[
1=\sum_{j=0}^N c_j^2\quad\textrm{and}\quad b_r=\sum_{j=0}^{N-r} c_j c_{j+r}\ \textrm{for}\ r=1,\ldots,N.
\]

When the number of non-zero elements of the sequence $\left(b_n\right)_n$ is infinite, from the
Szeg\"{o}--Kolmogorov--Krein theorem, under the additional assumption
\[
\int_{-\pi}^{\pi}\log b\left(\theta\right)  d\theta>-\infty,
\]
$b$ is non-negative if and only if
\[
b(\theta)=\left| C\left( \exp({i\theta})\right) \right|^{2},
\]
where
\[
C(z)=\sum_{n=0}^\infty c_n z^n=\exp\left(\frac{a_0}{2}+\sum_{n=1}^\infty a_n z^n\right),
\]
and
\[
a_n=\frac{1}{2\pi }\int_{-\pi }^{\pi }\exp \left( i n\theta \right) \log b\left( \theta \right) d\omega\quad
\textrm{for}\ n=0,1,2,\ldots .
\]
In this case, the relation between coefficients $\left(b_{n}\right)_n$ and $\left(c_{n}\right)_n$ is rather
complicated.

Nevertheless, it is easy to notice that if only $\sum_{n\in\nat}\left\vert b_{n}\right\vert<1/2$, then $b$ is positive.
(This condition, in the setting of \cite{brycannals}, reduces to $|a|<1/2$, see Remark \ref{brycannals}.)

\end{remark}

\section*{Acknowledgement}
The authors thank W.Bryc, V.Kaftal, M.Peligrad, R.Lata{\l}a and J.Weso{\l}owski for useful comments. The conversations
with W.Bryc and J.Weso{\l}owski were especially helpful. We also would like to thank the Referee for constructive
remarks that improved the paper.

The work on the revised version of the paper was done while the first named author (WM) was visiting the University of
Cincinnati.

\bibliographystyle{amsplain}
\bibliography{g:/moje/research/toolbox/mojebiby}

\end{document}